\documentclass[12pt,a4paper]{article}
\usepackage{amsmath, amsfonts, amssymb,amsthm}

\begin{document}
\author
{I.Kh. MUSIN}
\title {
ON THE FOURIER-LAPLACE TRANSFORM OF FUNCTIONALS
ON A WEIGHTED SPACE
OF INFINITELY DIFFERENTIABLE  FUNCTIONS}
\date{}
\maketitle
\newtheorem{theor}{Theorem}
\newtheorem{lemma}{Lemma}
\newtheorem*{sled}{Corollary}
\newtheorem{sleds}{Corollary}
\newtheorem{note}{Note}
\newtheorem*{theo}{Theorem S}
\newcommand{\eps}{\varepsilon}
\begin{abstract}
The strong dual space of linear continuous functionals on a weighted
space $ G $ of infinitely differentiable functions
defined on the real line is described in terms of their
Fourier-Laplace  transforms.
\end{abstract}

\vspace{0.3cm}

{\bf 1. Preliminaries and the main result.}

\vspace{0.3cm}

Let  $ M_0=1, M_1, M_2, \ldots,$ be an increasing
sequence of positive  numbers which satisfies the following conditions:

{\it i}).  $ M_k^2 \le M_{k-1} M_{k+1}  \ \ \forall k \in {\mathbb N}$;

{\it ii}).   $ \lim \limits_{k \to \infty}
{\left
(\displaystyle \frac {M_{k+1}} {M_k} \right)}
^{\frac 1 k} = 1.
$

{\it iii}).  $ \exists   Q_1, Q_2 > 0 $ such that
$$
M_k \ge Q_1 Q_2^k k! \ \ \forall k \in {\mathbb Z_+}.
$$

Let
$$
w(r)=
\sup_{k \in {\mathbb Z}_+}
\ln \frac {r^k}{M_k} \ , \  r > 0, \ w(0)=0.
$$
$w$ is continuous for $ r \ge 0 $ [3]
and $w(r) = 0 $ for  $ r \in [0, M_1] $. From this and
the condition {\it iii}) it follows that there exists
$A_w >0 $ such that
\begin{equation}
w(r) \le A_w r \ , \ r \ge 0.
\end{equation}
It is clear that  $w(\vert z \vert)$ is a subharmonic function in
the complex plane.

For $ z \in {\mathbb C}, \ t>0 $ denote by $ D(z,t) $ an open disk
of radius  $t$ about a point $z$ and by
$\partial D(z,t)$ its boundary.

Let  $ \alpha >1 $ and
$\psi:{\mathbb R} \to [0,\infty) $
be   a convex function satisfying the conditions:

1.  $ \exists A_{\psi}>0 $ such that
for any $x_1, x_2 \in {\mathbb R} $
\begin{equation}
\vert \psi (x_1) - \psi(x_2) \vert \le A_{\psi} (1+ \vert x_1 \vert + \vert
x_2 \vert)^{\alpha -1} \vert x_1
- x_2 \vert;
\end{equation}

2. $ \lim \limits_{
x \to \infty}
\displaystyle \frac {\psi (x)} {\vert x \vert} = +\infty $;

3. for any  $ z \in {\mathbb C}, \vert z \vert >1, $
$$
\mu_{\psi}
(D(z,t))
\le c_{\psi} {\vert z \vert}^{\alpha -1} t, \ t \in (0, \vert z \vert),
$$
where
$\mu_{\psi}$  is a measure associated by Riesz with subharmonic function
$\psi( Im \ z ) $, $c_{\psi}>0 $  is some constant.

For  a function $g: {\mathbb R} \to [0, \infty) $
such that
$ \displaystyle \frac {g(x)} {\vert x \vert} \to {+\infty} $
as $ x  \to \infty $ the Young transform   $g^*$ of $g$ is
defined by
$$
g^*(x) = \sup \limits_{y \in {\mathbb
R}}(xy - g(y)), \ x \in {\mathbb R}.
$$
$g^* $ is convex and
$ \displaystyle \frac {g^*(x)} {\vert x \vert} \to {+\infty} $
as  $  x \to \infty $ [10]. Moreover, if $ g $ is convex
then the inversion formula holds [10]:
 $(g^*)^* = g $.

Let  $\varphi  =  {\psi}^* $.
>From (2)
\begin{equation}
\varphi (x) > A_{\varphi} {\vert x \vert}^{\frac
{\alpha} {\alpha -1}} - B_{\varphi}, \ x \in {\mathbb R},
\end{equation}
where $A_{\varphi}, B_{\varphi}$ are some positive  numbers.

As usual ${\cal E} ({\mathbb R}) $ will
denote the  space of infinitely differentiable functions defined in
${\mathbb R}$ with the topology of uniform convergence of functions
and all their derivatives on compact subsets of   ${\mathbb R}$.
$ C({\mathbb R}) $ will denote   the space of continuous functions
in $ {\mathbb R} $.
$ H({\mathbb C})$ is the space of entire functions equipped with the
topology of uniform convergence on compact sets of complex plane.
For a locally convex space $X$ let
$X^*$  be  the strong dual of $X$.

Fix $\sigma >0 $. Let
$ \{{\varepsilon}_m \}_{m=1}^{\infty}$ be  an arbitrary
decreasing to zero sequence of positive  numbers.
 Set   $ {\theta}_m(x) =
 \exp (\varphi ( x ) - m \ln(1+\vert x \vert)),
x \in {\mathbb R}, m \in {\mathbb N} $. Let
$$
G_m =\{f \in {\cal E}({\mathbb R}): {\Vert f \Vert}_{G,m} =
 \sup_{x \in {\mathbb R}, k  \in {\mathbb Z_+}}
\displaystyle \frac {\vert f^{(k)}(x) \vert}
{(\sigma + {\varepsilon}_m)^k M_k
{\theta}_m(x)} < \infty \}, m \in {\mathbb N}.
$$
We let
 $G = \bigcap \limits_{m=1}^{\infty}  G_m
$
and endow this vector space with its natural projective limit topology.
It is clear that the definition of $G$ doesn't depend on the choice of
the sequence $\{\varepsilon_m \}_{m=1}^{\infty}$.
It is easy to show that for any $m \in {\mathbb N}$
the  canonical (inclusion) mapping
$ i_{m+1,m}: G_{m+1} \to G_m $ is relatively compact in $G_m$.
So $G$ is  the space $(M^*)$ (see for definitions [11])
>From this and  theorem 5 of [11] it follows that
  $G^*$ is a space  $(LN^*)$.

Let
$w_m(\vert z \vert)
= w ((\sigma +{\varepsilon}_m)^{-1} \vert z \vert),
z \in {\mathbb C}, m \in {\mathbb N}$.
Let
$$
P_m = \left\{f \in H({\mathbb C}):{\Vert f \Vert}_m =
\sup_{z \in  {\mathbb C}}
\displaystyle \frac {\vert f(z) \vert} {\exp(\psi( Im \ z ) +
w_m(\vert z \vert))} < \infty \right\}, m \in {\mathbb N}.
$$
Let  $P$ be the union of these normed spaces with a  topology of inductive
limit of the spaces  $P_m$.

For $T \in G^*$ we  define the
Fourier-Laplace transform
$ {\hat T }$ of $T$ by
$$
 {\hat T}(z) = T(e^{-ixz}) , \ z \in {\mathbb C}.
$$
The main point of our work is to characterize the dual space  $G^*$
of linear continuous functionals on $G$ in terms of their Fourier-Laplace
transforms.  The similar problem was considered in [14], but
the method presented  there is not applicable in  our case
(see condition $ \Phi $ 3 in [14], p. 48). Our principal result
is the following theorem.

\begin{theor} The Fourier-Laplace transform establishes topological isomorphism of
the spaces  $G^*$  and $P$.
\end{theor}
We establish  surjectivity  of the Fourier-Laplace transform
by using  the representation of  entire functions from $ P $ in the
form of the Lagrange series. This idea comes from [15].
To realise this idea we construct some special entire function.
For the  construction of this entire function we use results
of R.S. Yulmukhametov [16], [17] and M.I. Solomesch [7], [12], [13].

\vspace{0.3cm}

{\bf 2. Auxiliary results.}

\vspace{0.3cm}

\begin{lemma}
For any  $ z_1, z_2 \in {\mathbb C} $
\begin{equation}
\vert
w(\vert z_2 \vert) - w(\vert z_1 \vert) \vert \le A_w e
\vert z_2 - z_1 \vert.
\end{equation}
\end{lemma}
{\bf Proof}.
Let $ N(r)= \min \left\{k \in Z_+:
w(r) =  \ln \displaystyle \frac {r^k}{M_k} \right\}, \linebreak
r > 0 $.
Using the inequality from [15, Lemma 1.2]:
$
w(r) - w(1)
\ge N\left(\displaystyle
\frac r e \right)  , \ r > 1,
$
and the estimate (1), we obtain
\begin{equation}
 N(r) \le  A_w er, \ r > e^{-1}.
\end{equation}

>From the equality (see [3], [15, Lemma 1.2])
$$
w(r) - w(1) =
\int \limits_{1}^{r} \displaystyle \frac {N(x)} {x} \ dx \ , \ r > 1,
$$
and the estimate (5) we derive that:

$$
w(r_2) - w(r_1)  \le A_we (r_2 - r_1), \   r_2 > r_1 > 1;
$$

$$
w(r_2) - w(r_1 ) = w(r_2) -w(1)  \le A_we(r_2 -1)
\le  A_we (r_2 - r_1),   r_2 >  1, r_1 \in [0, 1].
$$
Since  $w(r) = 0 $ for  $ r \in [0, M_1] $ then from
two last inequalities the statement of Lemma follows.

\begin{lemma} For any  $ A >
0, m \in {\mathbb N} \exists  Q >0 $
such that
$$
w_m(\vert z \vert) + A \ln(1+\vert z \vert) \le
w_{m+1}(\vert z \vert) + Q, \ \ z \in {\mathbb C}.
$$
\end{lemma}

{\bf Proof}.
Let  $ N(r)$ be  defined as in Lemma 1.
Then  \ $ \displaystyle \frac {r^{N(r)}}{M_{N(r)}} \ge \displaystyle
\frac {r^{N(r)+1}}{M_{N(r)+1}} $. Consequently,
$r \le \displaystyle \frac {M_{N(r) +1}} {M_{N(r)}}$. From this, using
  {\it ii}),  for any $ \delta> 0 $ we can find
$A_{\delta } >0 $ such that
$r \le A_{\delta} (1 + \delta)^{N(r)}$.
Thus,
$N(r) \ge \displaystyle \frac {\ln \displaystyle \frac {r} { A_{\delta }}}
{\ln(1 + \delta )}$.
Set \  $ r= \displaystyle \frac {\vert z \vert}
{\sigma + {\varepsilon}_m} \ , \ z \in {\mathbb C}, z \ne 0 $. Then
\begin{multline*}
w_{m+1}(\vert z \vert) - w_m(\vert z \vert) \ge
\ln \left(\left(\displaystyle
\frac {\vert z \vert} {\sigma + {\varepsilon}_{m+1}}\right)^{N(r)}
\displaystyle
\frac 1 {M_{N(r)}}\right) - \\
\ln \left(\left(\displaystyle
\frac {\vert z \vert} {\sigma + {\varepsilon}_m}
\right)^{N(r)} \displaystyle
\frac 1 {M_{N(r)}}\right)
= k_0 (r)
\ln \displaystyle \frac
{\sigma + {\varepsilon}_m}
{\sigma + {\varepsilon}_{m+1}}  \ge
\displaystyle \frac {\ln \displaystyle \frac {\vert z \vert}
{ A_{\delta}(\sigma + {\varepsilon}_m)} }
{\ln(1 + \delta )}
\ln \displaystyle \frac
{\sigma + {\varepsilon}_m}
{\sigma + {\varepsilon}_{m+1}}.
\end{multline*}
Now we can  choose so small
$\delta $ and then  siutable constant  $Q \ge 0 $ so
that for all $ z \in {\mathbb C}$
$$
w_{m+1}(\vert z \vert) - w_m(\vert z \vert) \ge  A\ln(1+\vert z \vert)- Q.
$$

Lemma 2 is proved.

Using Lemma 2 it is easy to show that for any $m \in {\mathbb N}$
the  canonical (inclusion) mapping
$\gamma_{m, m+1}: P_m \to P_{m+1} $ is relatively compact  in $P_{m+1}$.
So $P$ is a space $(LN^*)$ (see for definitions [11]).

\begin{lemma} Let $T \in G^* $,  say for some
$  m_0 \in {\mathbb N}, c>0$
$$
\vert T(f) \vert \le c {\Vert f \Vert}_{G,m_0}, \ f \in G.
$$

Then $T$ can be represented in the form
$$
T(f) =
\sum \limits_{k=0}^{\infty}
\displaystyle \frac
{1} { (\sigma + {\varepsilon}_{m_0})^k M_k }
\int_{\mathbb R}
\displaystyle \frac {f^{(k)}(x)}{\theta_{m_0}(x)} \ d {\mu}_k (x) ,
$$
where ${\mu}_k $  are
  complex bounded measures
in ${\mathbb R}$ such that for all
$ k \in {\mathbb Z_+}$
$
\int_{\mathbb R} d \vert {\mu}_k \vert  (x)  \le  c.
$
\end{lemma}
The proof of the Lemma  is standard [1].

For
$z \in {\mathbb C}$   we let
$$
f_z(x)= \exp(-izx), \ x \in {\mathbb R}.
$$

\begin{lemma} For every
$z \in {\mathbb C}$ \ $f_z \in G$.
\end{lemma}

{\bf Proof}. For every $z \in {\mathbb C}, m \in {\mathbb N}$ we have
$$
{\Vert f_z \Vert}_{G,m} =
\sup \limits_{ x \in {\mathbb R}, k \in {\mathbb Z_+}}
									\displaystyle \frac {\vert
(-iz)^k \exp(-izx) \vert}
{(\sigma + {\varepsilon}_m)^k M_k
\theta_m(x)}
$$
$$
=
\exp(w_m(\vert z \vert)) \ \exp
(\sup \limits_{x
\in {\mathbb R}}(x \ Im  \ z - \ln \theta_m(x))),
\ z \in {\mathbb C}.
$$

Let's obtain the upper estimate
of $ \sup \limits_{x
\in {\mathbb R}}(x \ Im  \ z -
\ln \theta_m(x))
$.
>From  (3) it follows that for some
$A_m(\varphi), B_m(\varphi) >0 $
$$
\varphi(x) - m \ln (1 +  \vert x \vert ) \ge
 A_m(\varphi) {\vert x \vert}^{\frac
{\alpha} {\alpha -1}} - B_m(\varphi), \ x \in {\mathbb R}.
$$
Hence there exists  $y_m >0$ such that
for all $z \in {\mathbb C },
\vert Im \ z \vert > y_m,  $ \
the supremum of $ x \  Im \ z  - \varphi(x ) +m \ln (1 + \vert x \vert) $
over $ {\mathbb R}$  is  attained at some  point
$x_m \in (-(2A_m^{-1}(\varphi) \vert Im \ z \vert)^{\alpha -1},
(2A_m^{-1}(\varphi) \vert Im \ z \vert )^{\alpha -1})$.
If \  $ \vert Im \ z \vert > y_m,  $ \
we have
\begin{multline*}
\sup \limits_{x \in {\mathbb R}}
(x  \ Im \ z - \varphi(x) +m \ln (1 +  \vert x \vert)) -
\sup \limits_{x  \in {\mathbb R}}(x \ Im \  z - \varphi(x)) \\
\le
m \ln(1+ \vert x_m \vert) \le \ln (1 +
(2A_m^{-1}(\varphi) \vert Im \ z \vert )^{\alpha -1}) \\
\le m(\alpha -1) \ln (1 + \vert Im \ z \vert)
+ m(\alpha -1) \ln (1 + 2 A_m^{-1}(\varphi)) + m \ln 2.
\end{multline*}
Using inversion formula for the Young transform
and choosing  appropriate constant   $ b_m > 0 $, we have
for all  $z \in {\mathbb C}, m \in {\mathbb N}$
\begin{equation}
\sup \limits_{x  \in {\mathbb R}}
(x \ Im \ z  - \ln \theta_m(x)) - \psi ( Im \ z) \le
m(\alpha -1) \ln (1 + \vert Im \ z \vert) + b_m.
\end{equation}
Consequently,
$$
{\Vert f_z \Vert}_{G,m} \le
\exp(\psi ( Im \  z  )
+
w_m(\vert z \vert )+
m(\alpha -1)
\ln (1 + \vert Im \ z \vert) + b_m), z \in
{\mathbb C}.
$$

By Lemma 2  $ \exists q_m \ge 0 $ such that
\begin{equation}
{\Vert f_z \Vert}_{G,m} \le \exp(\psi ( Im \  z  ) +
w_{m+1}(\vert z \vert ) + q_m),
z \in {\mathbb C}.
\end{equation}
Thus,   $f_z \in G$ for all
$z \in {\mathbb C}$.

\begin{lemma}
Suppose that for $T \in G^*$ there are
  $c>0, m \in {\mathbb N}$ such that
\begin{equation}
\vert T(f)\vert \le c {\Vert f \Vert}_{G,m} \ ,  \ f \in G.
\end{equation}

Then $ {\hat T}$ is  entire  and satisfies
\begin{equation}
\vert
{\hat T}(z)
\vert \le c \exp(\psi(Im \ z)  + w_{m+1} (\vert
z \vert )+ q_m), \ z \in {\mathbb C},
\end{equation}
where  $ q_m $  is the same as in (7).
Moreover, ${\hat T}$ can be represented in the form
$$
{\hat T}(z) =\sum \limits_{k=0}^{\infty} V_k(z) z^k, \ z \in {\mathbb C},
$$
where $V_k $ are entire functions such that
for some  $A>0$ independent of $ k$
\begin{equation}
\vert V_k(z) \vert \le  \displaystyle \frac {A
(1 + \vert z \vert )^{m (\alpha -1)} \exp(\psi(Im \ z))}
{ (\sigma + {\eps}_m)^k M_k }, \ z \in {\mathbb C}.
\end{equation}
\end{lemma}

{\bf Proof}.
Using lemma 3 we have
\begin{equation}
{\hat T }(z) =
\sum \limits_{k=0}^{\infty}
\displaystyle \frac
{(-iz)^k}
{(\sigma + {\eps}_{m})^k M_k }
\int_{\mathbb R}
\displaystyle \frac {e^{-izx} }{\theta_{m}(x)} \ d {\mu}_k (x),
\ z \in {\mathbb C},
\end{equation}
where   ${\mu}_k $  are complex bounded measures
in ${\mathbb R}$ such that for all $ k \in {\mathbb Z_+}$

\begin{equation}
\int_{\mathbb R} d \vert {\mu}_k \vert  (x)  \le  c.
\end{equation}

Using the boundedness of the measures  ${\mu}_k $ and the estimate
\begin{equation}
\left\vert \displaystyle \frac {e^{-izx} }{\theta_{n}(x)} \right\vert \le
\exp(\psi ( Im \  z ) + {l}(\alpha -1) \ln (1 + \vert Im \ z \vert) +
b_{n}) ,
\end{equation}
which follows from (6) and holds
for every $ z \in {\mathbb C}, x \in {\mathbb R}, n \in {\mathbb N}$ ,
it's easy to show that
$$
V_k(z)= \displaystyle \frac
{(-i)^k}
{ (\sigma + {\eps}_{m})^k M_k }
\int_{\mathbb R}
\displaystyle \frac {e^{-izx} }{\theta_{m}(x)}
\ d {\mu}_k (x), \ z \in {\mathbb C}, k \in {\mathbb Z}_+,
$$
is entire. The estimate  (10) immediately  follows from  (12), (13)
(set $n=m$). From this  and  the condition {\it iii}) \ it follows
that the series  in the right-hand side of  (11) converges
in the topology of $H({\mathbb C})$.
Consequently, $ {\hat T}$ is entire. The estimate (9) is obtained
from (7) and (8).

\begin{lemma}
Let  $F \in P$ has the form
$F(z) = U(z) V(z),
z \in {\mathbb C}$,
where entire functions   $U, V $ for some $ m \in {\mathbb N}$,
 $C_U, C_V >0 $ satisfy the  estimates:
$$
\vert U(z) \vert
\le
\displaystyle \frac
{C_U
\exp(\psi(Im \ z )) }{1 + {\vert z \vert}^2} \ , \  z \in {\mathbb C};
$$
$$
\vert V(z) \vert \le C_V
e^{w_m(\vert z \vert)}, \  z \in {\mathbb C}.
$$

Then there exists  $T \in G^*$ such that $\hat T = F$
and
$$
\vert  T(f) \vert \le {\beta}_m C_U C_V {\Vert f \Vert}_{G, m+1} \ , \ f\in G,
$$
where ${\beta}_m >0 $ is some constant depending only on  $m$.
\end{lemma}

{\bf Proof}.
>From the estimate of $ \vert U(z) \vert $
it follows (see, for example, [4], [5], [8]) that
$ \exists p \in C({\mathbb R}) $ such that:

1. \ $  \sup \limits_{t \in {\mathbb R}}
\vert p(t) \vert \exp (\varphi(t)) < \displaystyle \frac {C_U} 2 $;

2. \ $ U(z) = \int_{\mathbb R} p(t) e^{-izt} \ dt, \
z \in {\mathbb C}$.

Further, by the Cauchy inequality for  Taylor coefficients
of  entire function
$
V(z) = \sum \limits_{k=0}^{\infty} v_k z^k $  we have
$$
\vert v_k \vert \le C_V
\inf \limits_{r>0}
\displaystyle \frac {\exp(w_m(r))}{r^k}
=
C_V (\sigma + {\varepsilon}_m)^{-k} \inf \limits_{r>0}
\displaystyle \frac {\exp(w(r))}{r^k} \ , \ k \in {\mathbb Z}_+.
$$
Using the equality [3]
\begin{equation}
\inf \limits_{r>0} \displaystyle \frac {\exp (w(r))}
{r^k}=
\displaystyle \frac
1 {M_k}\ , \ k \in {\mathbb Z}_+,
\end{equation}
we get
\begin{equation}
\vert v_k \vert \le C_V
{( \sigma + {\varepsilon}_m) }^{-k} {M_k}^{-1}, \ k \in {\mathbb Z}_+.
\end{equation}

Define the functional $T$ on  $G$ by the formula
$$
T(f) =
\int_{\mathbb R} p(t)  \sum \limits_{k=0}^{\infty} v_k i^k f^{(k)}(t) \ dt,
\ f \in G.
$$
It is  defined correctly. Indeed, for every $f \in G$
$$
\vert
\sum \limits_{k=0}^{\infty} i^k
v^k f^{(k)}(t) \vert
\le
\sum \limits_{k=0}^{\infty} \vert v_k \vert
{\Vert f \Vert}_{G,m+1}
(\sigma + {\varepsilon}_{m+1})
^k M_k
{\theta}_{m+1}(t), \ t \in {\mathbb R}.
$$

Using (15), we have for every $ f \in G $
$$
\vert \sum \limits_{k=0}^{\infty} v_k i^k
f^{(k)}(t) \vert \le {\beta}_m C_V
{\Vert f \Vert}_{G,m+1}{\theta}_{m+1}(t), \ t \in {\mathbb R},
$$
where  ${\beta}_m =
\displaystyle \frac {\sigma + {\varepsilon}_m}
{{\varepsilon}_m - {\varepsilon}_{m+1}}$ \ .
Consequently,
$$
\vert T(f) \vert \le {\beta}_m C_V
{\Vert f \Vert}_{G,m+1}
\int_{\mathbb R}
\vert p(t) \vert \exp(\varphi(t))(1+ \vert t \vert )^{-(m+1)} \ dt
$$
$$
\le {\beta}_m C_U C_V
{\Vert f \Vert}_{G,m+1}, \ f \in G.
$$
Thus, $T \in G^*$. Obviously,
$$
{\hat T}(z) = \int_{\mathbb R} p(t)
\sum \limits_{k=0}^{\infty}  v_k i^k (-iz )^k e^{-izt} \ dt
=
U(z) V(z)= F(z), \ z \in {\mathbb C}.
$$

\vspace{0.3cm}

 {\bf 3. The space  $ {\cal E}(\Phi)$.}

\vspace{0.3cm}

Let
$$
{\cal E}(\Phi)
= \{f \in {\cal E}({\mathbb R}):
\forall n \in {\mathbb Z_+}, m \in {\mathbb N} \ \ {\Vert f \Vert}_{n,m} =
\sup_{x \in {\mathbb R}, 0 \le k \le n}
\displaystyle \frac
{\vert f^{(k)}(x) \vert}
{{\theta}_m(x)} <  \infty \}.
$$
It is easy to see that  $ {\cal E}^*(\Phi) \subset G^*$.
In this section we study the
Fourier-Laplace transform of functionals from $ {\cal E}^*(\Phi)$.

If  $ F \in  {\cal E}^*(\Phi)$ then
$\exists $ $ c>0, n, m \in {\mathbb N} $ such that
$ \forall f \in {\cal E}(\Phi)$ \   we have
$ \vert F(f) \vert \le c {\Vert f \Vert}_{n,m} $.
So, using  (15), we have for some $C>0$
$$
\vert F(\exp(-ixz)) \vert \le C  (1 + \vert z \vert)^{n + m(\alpha -1)}
\exp (\psi(Im z)), \ z \in {\mathbb C}.
$$
By Lemma  5 \ $\hat F $ is an entire function.
The description of functionals from $ {\cal E}^*(\Phi)$
in terms of their Fourier-Laplace transforms will be done
in theorem 2. But at first we study the density of  polynomials
in  ${\cal E}(\Phi)$. We need the following simple lemma.

\begin{lemma} Let function  $g$ defined on the real line satisfies
for some constants  $A >0, B$ the inequality
$$
g(x) >
 A{\vert x \vert}^{\frac{\alpha} {\alpha -1}} - B
, \ x \in {\mathbb R}.
$$
Then there exists constant $C$ depending only on  $A$ and  $B$
such that for any   $ m \in {\mathbb  N}$
$$
\sup \limits_{x \in {\mathbb R}}
(m\ln (1+ \vert x \vert)- g(x))
<
(1-{\alpha}^{-1})m\ln m
+
Cm.
$$
\end{lemma}
{\bf Proof}. For any   $ m \in {\mathbb  N}$
$$
\sup \limits_
{x \in {\mathbb R}}
(m \ln (1+ \vert x \vert)
- g(x))
\le
\sup \limits_
{x \in {\mathbb R}}
(m \ln (1+ \vert x \vert)
-  A{\vert x \vert}^{\frac{\alpha} {\alpha -1}}) + B
$$
$$
=
B +
\sup \limits_
{x \ge 0}
(m \ln (1+  x )
-  A x^{\frac{\alpha} {\alpha -1}})
<
B+ 2 m +
\sup \limits_{u > 0}
(m u -  A e^{ {\frac{\alpha} {\alpha -1}}u})
$$
$$
= B+ 2 m + (1-{\alpha}^{-1})m\ln m
-
(1-{\alpha}^{-1})
m\ln ((1-{\alpha}^{-1}) A e) .
$$

\begin{theor}
The polynomials are dense in $ {\cal E}(\Phi) $.
\end{theor}

{\bf Proof}.
Let  $ f \in {\cal E}(\Phi) $, that is  $f \in
{\cal E}({\mathbb R})
$ and
for any
$n, m \in {\mathbb N}$
there exists  $c_{m,n} >0$ such that
for all
$
x \in {\mathbb R}
, k =0, 1, \ldots , n
$
\begin{equation}
\vert f^{(k)}(x) \vert \le c_{m,n} \theta_m(x).
\end{equation}
	Let us approximate $ f $ by polynomials in ${\cal E}(\Phi)$.
There are three steps in the proof.

1. Let   $\gamma \in {\cal E}({\mathbb R})$  be such that
supp   $ \gamma \subseteq [-2,2],
\gamma (x) =1 $ for  $ x \in [-1,1]$, $ 0 \le  \gamma (x) \le 1 $
$ \forall x \in {\mathbb R}$. Set
$ f_{\nu}(x) = f(x)\gamma (\frac x {\nu}), \ {\nu}
\in {\mathbb N},
x \in {\mathbb R} $. Obviously, $f_{\nu} \in
{\cal E}(\Phi)$.

Let us show that
$f_{\nu}$
 $\to $ $f$ in  $
{\cal E}(\Phi)
$ as  ${\nu} \to \infty$.
Take arbitrary  $m, n \in {\mathbb N}$. Then
$$
\sup \limits_{x \in {\mathbb R}} \displaystyle \frac
{\vert f_{\nu}(x) - f(x) \vert} {\theta_m(x)} \le
\sup \limits_{\vert x \vert > {\nu}}
\displaystyle \frac {\vert f(x)
\vert}{{\theta}_m(x)}
\le
\sup \limits_{\vert x \vert > {\nu}}
\displaystyle \frac {c_{m+1,n} {\theta}_{m+1}(x)}
{{\theta}_m(x)}.
$$

Consequently,  as
${\nu} \to \infty$
\begin{equation}
\sup \limits_{x \in {\mathbb R}} \displaystyle \frac
{\vert f_{\nu}(x) - f(x) \vert} {\theta_m(x)} \to 0.
\end{equation}
Then,
$$
\sup \limits_{x \in {\mathbb R},
1 \le k \le n
}
\displaystyle \frac
{\vert (f_{\nu}(x)
 - f(x))^{(k)}
\vert} {\theta_m(x)}
$$
$$
= \sup \limits_
{x \in {\mathbb R},
1 \le k \le n
}
\displaystyle \frac {\vert \sum
\limits_{s=o}^{k-1} C_k^s f^{(s)}(x) {\nu}^{s-k} \gamma^{(k-s)}
( \frac x {\nu}) + f^{(k)}(x) (\gamma (\frac x {\nu}) -1) \vert}
{\theta_m(x)}
$$
$$
\le
\sup \limits_{ \nu < \vert x \vert < 2{\nu},
1 \le k \le n}
\displaystyle \frac {\sum
\limits_{s=o}^{k-1} C_k^s \vert f^{(s)}(x)\vert  n^{s-k}
\vert \gamma^{(k-s)}( \frac x {\nu})\vert}
{\theta_m(x)}
+ \sup \limits_{\vert x \vert >
{\nu},
1 \le k \le n}
\displaystyle \frac {\vert f^{(k)}(x) \vert}
{\theta_m(x)}.
$$
>From this, using (16) (substititing  $m$ by $m+1$),
 we conclude that as  $ {\nu} \to \infty$
$$
\sup \limits_{x \in {\mathbb R},
1 \le k \le n
}
\displaystyle \frac
{\vert (f_{\nu}(x)
 - f(x))^{(k)}
\vert} {\theta_m(x)}
\to 0.
$$
>From this and (17) it follows that
$ {\Vert f_{\nu} - f \Vert}_{n,m} \to 0 $
as  $ \nu \to \infty $.
Since  $ m, n \in {\mathbb N}$  are arbitrary
this means that
the sequence
 $ \{f_{\nu}\}_{\nu=1}^{\infty} $ converges
to  $f$ in
${\cal E}(\Phi)$
as  $ \nu \to \infty $.

2. Fix   $ \nu \in {\mathbb N}$.
Let
$ h(z) = \displaystyle \sum \limits_{k=0}^{+\infty} a_k z^{2k},
z \in {\mathbb C},  h \not \equiv 0, $  be an entire function of
 exponential type 1 such that
$ h \in L_1({\mathbb R}),  h(x) \ge 0, x \in  {\mathbb R} $. For example,
we may put  $ h(z) =
\displaystyle \frac  {\sin^2{ \frac z 2}} {z ^2}, \
z \in~{\mathbb C} $.

By the Paley-Wiener theorem
$ \exists g \in C({\mathbb R})$ with   support  $ g \subseteq [-1,1] $
such that
$$
h(z) = \int \limits_{-1}^1 g(t) e^{-izt} \ d  t \ , \ z \in {\mathbb C}.
$$

Since
$$
h^{(k)}(z) = \int \limits_{-1}^1 g(t) (-it)^k e^{-izt} \ d t \ ,
\ z \in {\mathbb C},
\ k \in {\mathbb Z_+},
$$
then
\begin{equation}
\vert h^{(k)}(x) \vert \le C_g \ , \ x \in {\mathbb R},
\ k \in {\mathbb Z_+},
\end{equation}
where
$C_g= 2\max \limits_{\vert t \vert \le 1} \vert g(t) \vert $.

Let
$ \int \limits_{-\infty}^{+\infty} h(x) \ d x = A $. For  $\lambda>0$ we set
$$
f_{\nu,\lambda}(x) =  \displaystyle \frac  {\lambda} A \int_{\mathbb R}
f_{\nu}(y)
h(\lambda (x-y)) \ d y, \ x \in {\mathbb R}.
$$
It is easy to see that
$
f_{\nu,\lambda}
\in
{\cal E}(\Phi)$.

Let us show that
$
f_{\nu,\lambda}
\to
f_{\nu}
$ in  $
{\cal E}(\Phi)$
as
$ \lambda \to +\infty $. Take arbitrary  $m, n \in {\mathbb N}$.
For any  $k  \in {\mathbb Z}_+, x \in {\mathbb R}$
$$
f_{\nu,\lambda}
^{(k)}
(x) -
f_{\nu}
^{(k)}(x) =
\displaystyle \frac  {\lambda} A
\int_{\mathbb R} (
f_{\nu}
^{(k)}(y)-
f_{\nu}
^{(k)}(x))h(\lambda (x-y)) \ dy
$$
$$
= \displaystyle \frac  {\lambda} A
\int_{\vert y -x \vert \le {\lambda}^{-\frac 2 3}}
(
f_{\nu}
^{(k)}(y)-
f_{\nu}
^{(k)}(x))h(\lambda (x-y)) \ d y
$$
$$
+ \displaystyle \frac  {\lambda} A
\int_{\vert y -x \vert > {\lambda}^{-\frac 2 3}}
(
f_{\nu}
^{(k)}(y)-
f_{\nu}
^{(k)}(x))h(\lambda (x-y)) \ d y =
I_{1,k}(x) + I_{2,k}(x).
$$
Let $ K_n= \max
\limits_{x \in {\mathbb R}, 0 \le k \le n+1}
\vert f_{\nu}^{(k)}(x) \vert $.
Obviously,
\begin{equation}
\vert I_{1,k}(x) \vert
\le \displaystyle \frac
{2C_gK_n}{ A}
{\lambda}^{- \frac 1 3}
\ , \ x \in {\mathbb R}, \ k =0, 1, \ldots , n;
\end{equation}
\begin{equation}
\vert I_{2,k}(x) \vert \le
\displaystyle \frac
{2K_n}
 A
\int_{\vert t \vert >
{\lambda}
^{\frac 1 3}}
h(t) \ d t \  , \ x \in {\mathbb R},
\ k =0, 1, \ldots , n.
\end{equation}
Let  $\varepsilon > 0 $ be arbitrary and
$\Theta_m = \inf \limits_
{ x \in {\mathbb R}}
{\theta_m(x)}, m~{ \in {\mathbb N}} $.
Choose
$
{\lambda}
(\varepsilon )
>0$ so that
$
\displaystyle \frac
1 A
\int_{\vert t \vert >
{\lambda}
^{\frac 1 3}}
h(t)  \ d t <
\displaystyle \frac
{\varepsilon \Theta } {4K_n }
$ \ and
\ $
{2C_gK_n}
{\lambda}^{- \frac 1 3} <
\displaystyle
\frac {A\varepsilon \Theta } {2} $
for  ${\lambda} > {\lambda}(\varepsilon)$.
Then from (19) and (20) it follows that
$ {\Vert f_{\nu,\lambda} - f_{\nu} \Vert}_{n,m} < \varepsilon $
for all ${\lambda} > {\lambda}(\varepsilon )$.
Consequently, by the semi-norm
${\Vert \cdot \Vert}_{n,m}$ \
$f_{\nu,\lambda}\to f_{\nu} $ as  $ \lambda \to +\infty $.
Since  $n, m \in {\mathbb N} $ are arbitrary, it means   that
$
f_{\nu,\lambda}
\to f_{\nu}
$ in  $
{\cal E}(\Phi)
$ as  $ \lambda \to +\infty $.

3. Fix $\lambda>0, {\nu} \in {\mathbb N} $.  Let us approximate
$
f_{\nu,\lambda}
$
by the polynomials
in ${\cal E}(\Phi)$.

Let  $h(x) =
\displaystyle \sum \limits_{k=0}^{+\infty} a_k x^{2k},
x \in {\mathbb R}$.
Let  $P_{2N}(x) =
\displaystyle \sum \limits_{k=0}^{N} a_k x^{2k}, x \in {\mathbb R},
N  \in {\mathbb N}$.
Using the Taylor formula and the estimate (18), we have
\begin{equation}
\vert h(x) - P_{2N}(x) \vert \le
\displaystyle \frac {C_g {\vert x \vert}^{2N+1}}{(2N+1)!}, \ x \in {\mathbb R}.
\end{equation}

Set
$$
Q_{2N}(x) =
\displaystyle \frac  {\lambda} A
\int
\limits_{-L}^L
f_{\nu}(y) P_{2N}(\lambda (x-y)) \ d y , \ x \in {\mathbb R},
N
\in {\mathbb N},
$$
where $L>0$ is chosen so that   supp $ \
f_{\nu}
\subseteq [-L,L]$.
$Q_{2N}$ is a polynomial of the degree  not more than $2N$.

Let us show that the sequence
$\{Q_{2N}\}_{N=1}^{\infty} $ converges in  $
{\cal E}(\Phi)
$ to $
f_{\nu,\lambda}
$  as
$N \to \infty $.
We take   arbitrary $m, n \in {\mathbb N}$.
For any $k \in {\mathbb Z}_+ $
$$
f_{\nu,\lambda}
^{(k)}(x) - Q_{2N}^{(k)}(x)
=  \displaystyle \frac {\lambda} A
\int \limits_{-L}^L
f_{\nu}
^{(k)}(y)
(h({\lambda}(x-y)) -
P_{2N}(\lambda (x-y))) \ d y , \ x \in {\mathbb R}.
$$
Using (18), (21) we have for all
$ x \in {\mathbb R}, N {\in {\mathbb  N}}, k=0, 1, \ldots , n $

$$
\vert
f_{\nu,\lambda}
^{(k)}(x) - Q_{2N}^{(k)}(x) \vert \le K_n C_g
\displaystyle \frac
{2L {\lambda}^{2N+2}} A
 \displaystyle \frac {(\vert x \vert + L)^{2N+1}} {(2N+1)!} \ .
$$
Thus, there exist   $C_1>0 $ depending on  $n$ and  $C_2 $
such that for any
$ N {\in {\mathbb  N}}$
\begin{equation}
{\Vert
f_{\nu,\lambda}
- Q_{2N} \Vert}_{n,m} \le
\displaystyle \frac
{C_1 C_2^{2N+2}}{(2N+1)!}
\sup \limits_{x \in {\mathbb R}}
\displaystyle \frac {(1+\vert x \vert)^{2N+1}} {\theta_m(x)} \ .
\end{equation}

Since $ \varphi $  satisfies the conditions of Lemma 7 then
there exists  positive number
$C_3 $
such that for all  $ m, N \in {\mathbb  N} $
$$
\sup \limits_{x \in {\mathbb R}}
((2N+m+1)\ln (1+ \vert x \vert)-
\varphi (x))
$$
$$
\le (1- {\alpha}^{-1})(2N+m+1) \ln(2N+m+1) + C_3 (2N+m+1).
$$
>From this estimate and the inequality (22) it follows that
there exist positive numbers  $C_4, C_5$ depending on m, n
 such that  for all $  N \in {\mathbb  N} $
$$
{\Vert
f_{\nu,\lambda}
- Q_{2N} \Vert}_{n,m} \le
 \displaystyle \frac {C_4 C_5^N (2N+1)^{(2N+1)(1 - {\alpha}^{-1})}}
{(2N+1)!} \ .
$$
The right-hand side of the last inequality
tends to  0 as $N \to \infty $.
This means that
 $
f_{\nu,\lambda}
$ is approximated
by polynomials in  $
{\cal E}(\Phi)
$
since $m, n \in {\mathbb N} $ are arbitrary.

The density of polynomials in ${\cal E}(\Phi)$ follows
from  steps 1) -- 3).

\begin{theor}
Let   $  U \in H({\mathbb C}) $
for some  $ C>0, N \in {\mathbb N} $  satisfies the eniquality
$$
\vert U(z) \vert \le C
( 1 + \vert z \vert)^N \exp (\psi(Im z)), \ z \in {\mathbb C}.
$$
Then there exists the unique  functional
$T \in
{\cal E}^*(\Phi)$
such that  $ T(\exp(-izx)) = U(z), z \in {\mathbb C} $ and
$$
 \vert T(f) \vert \le A C {\Vert f \Vert}_{N+2,2} \ , \
f \in{\cal E}(\Phi),
$$
where  $A >0 $  doesn't depend on  $ U $.
\end{theor}

{\bf Proof}.
Let us take real numbers
$ {\lambda}_0, {\lambda}_1, \ldots , {\lambda}_{N+1} $
not equal to each other and zero.
Choose numbers  $a_0, a_1, \ldots, a_{N+1}$ so that an  entire function
$$
g(z) = U(z) -
\sum \limits_{k=0}^{N+1} a_k \exp(-i{\lambda}_k z),
\  z \in {\mathbb C},
$$
satisfies the following condition:
 $ g^{(n)}(0) = 0, n = 0, 1, \ldots , N+1. $
It is clear that these numbers are uniquely defined.
>From the estimate on   $U$ we have
 \ $\vert U^{(k)}(0) \vert \le c_1 C, k =0, 1, \ldots , N+1, $
where  $c_1>0$ depends only on  $N$.
Then
$\vert a_k \vert \le c_2 C,
\ k =0, 1, \ldots , N+1, $
where  $c_2>0$ depends only on  $N$.
Using the second condition on  $\psi$ we have for some  $c_3 >0$
depending only on  $N$
$$
\vert g(z) \vert \le
c_3 C (1 + \vert z \vert)^N \exp (\psi(Im z)), \  z \in {\mathbb C}.
$$
Let
$h(z) =
\displaystyle \frac
{g(z)} {(-iz)^{N+2}}$.
Obviously,
$$
\vert h(z) \vert \le
\displaystyle \frac {c_4 C\exp (\psi(Im z))}
{1 + {\vert z \vert}^2} \ ,
\  z \in {\mathbb C},
$$
where  $c_4>0$  depends only on  $N$.

As in Lemma 6  $ \exists p \in C({\mathbb R})$ such that:

1. \ $  \sup \limits_{t \in {\mathbb R}}
\vert p(t) \vert \exp (\varphi(t)) < \displaystyle \frac {c_4 C} 2 $;

2. \ $ h(z) =
\int_{\mathbb R} p(t) e^{-izt} \ dt, \
z \in {\mathbb C}$.

Define functional  $T$ on
$
{\cal E}(\Phi)
$
by the formula
$$
T(f) =
\sum \limits_{k=0}^{N+1}
a_k f({\lambda}_k) +
\int_{\mathbb R} p(t) f^{(N+2)}(t) \ dt, \ f \in {\cal E}(\Phi).
$$
It is easy to show that for some   $A >0 $ depending only on  $N$
$
 \vert T(f) \vert \le A C {\Vert f \Vert}_{N+2,2} \ , \
f \in{\cal E}(\Phi).
$
Consequently,
$ T \in
{\cal E}^*(\Phi)$. Obviously,
$ T(\exp(-izx)) = U(z), z \in {\mathbb C} $.

Now we prove the uniqueness. Suppose that
$S \in{\cal E}^*(\Phi)$ \
$ S(\exp(-izx)) = 0 \ \forall z \in {\mathbb C} $.
Let us show that   $ S = 0 $.
Let $ P_{N,\xi}(x) = \sum \limits_{\nu=0}^{N}
\displaystyle \frac {{\xi}^{\nu}} {{\nu}!} x^{\nu},
\xi,  x \in {\mathbb R}, \ N \in {\mathbb Z_+}$.
We show at first that for any $R>0$ the sequence of polynomials
$ \{P_{N,\xi}\}_{N=0}^{\infty}$ converges in ${\cal E}(\Phi)$  to
$e^{\xi x}$ uniformly by $\xi \in [-R,R]$ as  $N \to \infty $.
Let   $n, m \in {\mathbb N} $ be arbitrary.
>From the Taylor formula
$$
\vert
e^{\xi x} -  P_{N-k,\xi}(x)
\vert \le
\displaystyle \frac
{{\vert \xi \vert}^{N-k+1}
{\vert x \vert }^{N-k+1}}
{
(N-k+1)!}
 \max(1, e^{\xi x})
, \ k= 0, 1, \ldots ,  N,
$$
then as  $N \ge n $
$$
\sup_{x \in {\mathbb R}, 0 \le k \le n}
\displaystyle \frac
{\vert (e^
{\xi x})^{(k)}
-  P_{N,\xi}^{(k)} (x) \vert}
{\theta_m(x)}
=
\sup_{x \in {\mathbb R}, 0 \le k \le n}
\displaystyle \frac
{{\xi}^k \vert e^{\xi x} -  P_{N-k,\xi}(x)\vert}
{\theta_m(x)}
$$
$$
\le
\displaystyle \frac{{\vert \xi \vert}^{N+1}}
{(N-k+1)!}
\sup \limits_{x \in {\mathbb R}}
\displaystyle \frac
{\vert x \vert^{N- k +1} \max(1, e^{\xi x})}
{\theta_m(x)}
$$
$$
\le
\displaystyle \frac{R^{N+1}}
{(N-n+1)!}
\sup \limits_{x \in {\mathbb R}}
\displaystyle \frac
{(1 + \vert x \vert)^{N+1}
\max(1, e^{\xi x})}
{\theta_m(x)}.
$$
Since $ \exists A, B >0 $ depending  on  $R$ such that for all  $\xi
\in [-R,R], \ x \in {\mathbb R},
\
\varphi(x) - \xi x
\ge A {\vert x \vert}^{\frac
{\alpha} {\alpha -1}} - B,
$
then, using Lemma  7, we can find
constant  $Q>1$ independent of $ \xi \in [-R,R]$ such that
$$
\sup \limits_{x \in {\mathbb R}}
\displaystyle \frac
{(1 + \vert x \vert)^{N+1}
\max(1, e^{\xi x})}
{\theta_m(x)}
\le
Q^
{N+m+1}
(N+m+1)
^{(1-{\alpha}^{-1})(N+m+1)}.
$$
>From these two last estimates it follows that
$$
\sup_{x \in {\mathbb R}, 0 \le k \le n}
\displaystyle \frac
{\vert (e^
{\xi x})^{(ek)}
-  P_{N,\xi}^{(k)} (x) \vert}
{\theta_m(x)}
\to 0
$$
uniformly by  $\xi \in [-R,R]$ as  $ N \to \infty $.
Since  $n, m $ are arbitrary, then it means that for any $R>0$
$ \{P_{N,\xi}\}_{N=0}^{\infty} \to e^{\xi x}$ in
${\cal E}(\Phi)$  uniformly by  $\xi \in [-R,R]$ as $N \to \infty $.
>From this we conclude  that
$
\sum \limits_{\nu=0}^{\infty}
\displaystyle \frac
{S(x^{\nu})}
{{\nu}!}
{\xi}^{\nu}
$
converges uniformly on compacts from   ${\mathbb R}$ to $S(\exp(\xi x))$.
By assumption
for any $ \xi \in {\mathbb R}$ \ $ S(\exp(\xi x))= 0 $. Consequently,
$S(x^{\nu})= 0 $ for all  $\nu \in {\mathbb Z}_+$. Since the
polinomials are dense in  ${\cal E}(\Phi)$ then  $ S=0 $.

 The theorem is proved.

\vspace{0.3cm}

 {\bf 4. The proof of theorem 1.}

\vspace{0.3cm}

 In this section at first we  briefly  describe a special case of
M. I. Solomesch's general result [7], [12].

Let  $f$ be  an arbitrary  entire function with zeros
$\{{\lambda}_j \}_{j=1}^{\infty} \subset {\mathbb C}$
of corresponding multiplicities
$\{m_j \}_{j=1}^{\infty} $. The Weierstrass representation
of  $f$ has the form [2, p. 17]
$$
f(\lambda) =
g(\lambda)
\displaystyle \prod \limits
_{j=1}^{\infty}
\left(1 -
\frac
{\lambda}
{\lambda}_j \right)
^{m_j}
\exp \left(P_j\left(
\displaystyle
\frac {\lambda}{{\lambda}_j}\right)\right), \lambda
\in
{\mathbb C},
$$
where  $g $ is an entire function without zeros,  $P_j$ are special
polynomials, $ j \in {\mathbb N}$.

Let $ t = \{t_j\}_{j=1}^{\infty} $ be a sequence of complex numbers such
that
for any
$j \in {\mathbb N}$ \ ${\lambda}_j + t_j \ne 0$.
Put
$$
f_t(\lambda)=
g(\lambda) \prod \limits
_{j=1}^{\infty}
\left(1 -
\frac
{\lambda}
{{\lambda}_j + t_j} \right)
^{m_j}
\exp \left(P_j\left(
\displaystyle
\frac {\lambda}{{\lambda}_j}\right)\right),  \lambda \in {\mathbb C}.
$$
Now we can formulate the result of M.I. Solomesch
[7],[12, chapter 2, $ \S 3 $].

\begin{theo} Let
the disks $ D({\lambda}_j, r_j) $ be  pairwise disjoint,  \
$ {\lambda}_j + t_j   \in D({\lambda}_j, r_j) $,
$ r_j > 0, j \in {\mathbb N} $,
\ and        \
$
\displaystyle \sum \limits_{j=1}^{\infty}
\displaystyle
\frac
{m_j \vert t_j \vert} {r_j} < \infty $.

Then  $f_t$ is an  entire function and for
  $\lambda $ outside
$\bigcup \limits_{j=1}^{\infty} D({\lambda}_j, r_j)$
$$
\vert \ln f_t (\lambda) \vert -
\vert \ln f (\lambda) \vert \vert \le C,
$$
where  $C>0$ is some constant.
\end{theo}

{\bf Proof of theorem 1}.

Obviously, the map $J: T \to  {\hat T}, \  T \in G^*$,
is linear and by Lemma  5 acts from  $G^*$ to  {$P$}.

$ J $ is continuous. Indeed, if
$T \in G^*$, then for some $ m \in {\mathbb N}$ \ $T \in G_m^* $. Hence,
$$
\vert T(f) \vert \le
{\Vert T \Vert}_{G^*,m}
{\Vert f \Vert}_{G,m} \ , \ f \in G,
$$
where
$
{\Vert T \Vert}_{G^*,m}
$ is a norm of  $T$ in  $G_m^*$. By Lemma  5
$$
\vert {\hat T}(z) \vert \le
{\Vert T \Vert}_{G^*,m}
\exp(
\psi( Im \ z )
+ w_{m+1}(\vert z \vert )+ q_m), \ z \in {\mathbb C},
$$
where $q_m >0 $ is  some constant. Consequently,
$$
{\Vert {\hat T} \Vert}_{m+1} \le
e^{q_m} {\Vert T \Vert}_{G^*,m} \ ,  \ T \in G_m^*.
$$

This means that  $J$ is continuous.

At first we  prove that  $J$  is injective. That is, if
for  $ T \in G^* \ \ {\hat T} \equiv 0 $, then
$T(f)=0 $ \ $ \forall f \in G $.
We will follow the scheme of the paper  [14, see the proof of theorem 2.8].
By Lemma 3 \ $T$  is represented in the form
$$
T(f) =
\sum \limits_{k=0}^{\infty}
\displaystyle \frac
{1} {(\sigma + {\eps}_{m})^k M_k }
\int_{\mathbb R}
\displaystyle \frac {f^{(k)}(x)}{\theta_{m}(x)} \ d {\mu}_k (x) , \ f \in
G,
$$
where  ${\mu}_k$ are complex bounded measures in  ${\mathbb R}$
such that
 $ \forall k \in~{ {\mathbb Z_+}}$
$
\int_{\mathbb R} d \vert {\mu}_k \vert  (x)  \le
{\Vert T \Vert}_{G^*,m}.
$
Note that functional
$$
T_k(f)= \displaystyle \frac
{1} {(\sigma +  {\eps}_m)^k M_k }
\int_{\mathbb R}
\displaystyle \frac {f(x)}{\theta_{m}(x)} \ d {\mu}_k (x) ,
\ f \in
{\cal E}(\Phi),
$$
is in  ${\cal E}^*(\Phi)$.
By Lemma  5
$
{\hat T}(z) =
\sum \limits_{k=0}^{\infty}
V_k(z) z^k $,
where entire functions
$V_k(z) = (-i)^kT_k(\exp(-ixz))$
for some
$B_m>0$ independent of $k$ satisfy the estimate
\begin{equation}
\vert V_k(z) \vert \le
\displaystyle \frac
{B_m (1 + \vert z \vert )^{m(\alpha -1)} \exp(\psi(Im \ z))}
{(\sigma +  {\eps}_m)^k M_k } \ , \ z \in {\mathbb C}.
\end{equation}
Consider
$H(z,u)
=  \sum \limits_{k=0}^{\infty} V_k(z) u^k,
\ z, u  \in {\mathbb C} $.
It is easy to see that
$H(z,u) $ is entire. From (23) we have an estimate for $H(z,u)$
$$
\vert H(z,u) \vert \le
\displaystyle \frac {B_m
(\sigma + {\eps}_m)}
{{\eps}_m - {\eps}_{m+1}}
(1 + \vert z \vert )^{m(\alpha -1)} \exp(\psi(Im \ z) +
w_{m+1}(\vert u \vert)), \  z, u
\in {\mathbb C}.
$$
By our hypothesis $ H(z,z) = 0, \ z  \in {\mathbb C} $.
Hence,  $ H(z,u) = (z-u) S(z,u)$, where  $ S $ is entire.
Note  that $S$ satisfies the  estimate
\begin{equation}
\vert S(z,u) \vert \le A_S
(1 + \vert z \vert )^{m(\alpha -1)} \exp(\psi(Im \ z) +
w_{m+1}(\vert u \vert))
, \  z, u  \in {\mathbb C},
\end{equation}
where  $A_S >0 $ is some constant.  It is obvious when
$ \vert z - u \vert \ge 1$ and  for   $ \vert z - u \vert < 1 $
it can be obtained by applying the maximum principle and
by using the inequality  (4). Then we expand $S$ in a power series in $u$:
$ S(z,u)=
\sum \limits_{k=0}^{\infty}
S_k(z) u^k $.
By Cauchy's inequality for the coefficients in a power series expansion,
equality (14), inequality (24), we have
$$
\vert S_k(z) \vert \le
\displaystyle \frac {A_S} {(\sigma + {\eps}_{m+1})^k M_k}
(1 + \vert z \vert )^
{m(\alpha -1)}
\exp(\psi(Im \ z)) \ ,
k \in {\mathbb Z}_+,
z  \in {\mathbb C}.
$$
By theorem  2 there exist ${\Phi}_k \in {\cal E}^*(\Phi)$
such that  $ {\Phi}_k(\exp(-izx)) = S_k(z), z\in {\mathbb C} $
and \begin{equation}
 \vert {\Phi}_k(f) \vert \le
\displaystyle \frac {K} {(\sigma +{\eps}_{m+1})^k M_k}
{\Vert f \Vert}_{N+2,2} \ , \
f \in{\cal E}(\Phi),
\end{equation}
where   $ N =[m(\alpha -1)] +1 $,
constant  $K >0 $  doesn't depend on  $ k \in {\mathbb Z}_+$.

Since  $H$  is represented in the form
$$
H(z,u) = z S_0(z)  +   \sum \limits_{k=1}^{\infty}
(z S_k(z)- S_{k-1}(z)) u^k,
$$
then
  $ V_0(z) = z S_0(z), V_k(z) = zS_k(z) - S_{k-1}(z), k \in {\mathbb N} $.
By theorem 2 \
$T_0(f)=i{\Phi}_0(f'), \ (-i)^k T_k(f)= i{\Phi}_k(f') - {\Phi}_{k-1}(f),
f  \in{\cal E}(\Phi)$.
According to
{\it ii}) $ \forall  \delta> 0  \ \exists  A_{\delta } > 0:
\  \forall n \in {\mathbb N} \ \
M_{n+1} \le A_{\delta } (1 + \delta)^n M_n $.
Then from  (25) we have for $ f \in G $
$$
 \vert
{\Phi}_n(f^{(n+1)})
\vert
\le
\displaystyle \frac
{K} {(\sigma +{\eps}_{m+1})^n M_n}
\sup \limits_{x \in {\mathbb R}, 0 \le k \le N+2}
\displaystyle \frac
{\vert f^{(n+k+1)}(x) \vert}{\theta_2(x)}
$$
$$
\le
\displaystyle \frac
{K} {(\sigma +{\eps}_{m+1})^n M_n}
\sup \limits_{x \in {\mathbb R}, 0 \le k \le N+2}
\displaystyle \frac  {
{\Vert f \Vert}_{G,m+2}
(\sigma + {\eps}_{m+2})^{n+k+1}
M_{n+k+1}
{\theta}_{m+2}(x)}
{{\theta}_2(x)}
$$
$$
\le
K_{\delta}
{\Vert f \Vert}_{G,m+2}
\left(\displaystyle \frac
{(\sigma + {\eps}_{m+2})(1 + \delta)^{(k+1)}}
{\sigma + {\eps}_{m+1}}\right)^n,
$$
where
$
K_{\delta}= K \sup \limits_{ 0 \le k \le N+2}
\left((\sigma + {\eps}_{m+2})^{k+1}
A_{\delta }^{k+1}
(1 + \delta)^{\frac {(k+1)k} 2}\right).
$
Choose  $\delta $ so small that
$\displaystyle \frac
{(\sigma + {\eps}_{m+2})(1 + \delta)^{(k+1)}}
{\sigma + {\eps}_{m+1}} < 1$.
Then  $  {\Phi}_n(f^{(n+1)})  \to 0 $ as $ n \to \infty $.
Therefore, for any  $f \in G$
$$
T(f)=
\sum \limits_{k=0}^{\infty}
T_k(f^{(k)})= i{\Phi}_0(f') +
\sum \limits_{k=1}^{\infty} (i^{k+1}{\Phi}_k(f^{(k+1)})-
i^k{\Phi}_{k-1}(f^{(k)}))
$$
$$=
\lim \limits_{n \to \infty} i^{n+1}{\Phi}_n(f^{(n+1)}) = 0.
$$

Now let us prove that   $ J $ is surjective.
Let an
entire function $F \in P $ for some $ m \in {\mathbb N}, c>0 $
satisfies
\begin{equation}
\vert F(z) \vert \le c \exp(
\psi( Im \ z ) +w_m(\vert z \vert)),
\ z \in {\mathbb C}.
\end{equation}
We wish to show that  there is $ {\cal F} \in G^* $ with
$ {\hat {\cal F}} =F$.

Since functions  $ \psi(Im \ z ), w_{m+3} (\vert \ z \vert)
$ satisfy the conditions of the theorem 6 of [16] (see also [17, theorem 4],
then  there are entire functions  $L$ and  $Y$
such that:

(L1).
All the zeros  $ {\{\lambda}_j \}_1^{\infty} $ of $ L $ are simple
and the disks
$ D({\lambda}_j, d_1 {\vert {\lambda}_j \vert}^{1-\alpha}) $
are disjoint
for some $d_1 > 0 $.

(L2).
Outside the set  $
\bigcup \limits_{j=1}^{\infty}
D({\lambda}_j, d_1 {\vert {\lambda}_j \vert}^{1-\alpha})
$
\begin{equation}
\vert
\psi(Im \ z )
  - \ln \vert L(z) \vert \vert
\le A \ln(1 + \vert z \vert ) + A_0,
\end{equation}
where $A, A_0$ are some positive  numbers.

(Y1). All the zeros  $ {\{\mu}_k \}_1^{\infty} $ of $ Y $ are simple
and the disks $ D({\mu}_k, d_2  ) $ are disjoint for some $ d_2 > 0 $.

(Y2).
Outside the set
$ \bigcup \limits_{k=1}^{\infty} D({\mu}_k, d_2 )$
\begin{equation}
\vert
w_{m+3} (\vert \ z \vert)
  - \ln \vert Y(z) \vert \vert
\le B \ln(1 + \vert z \vert ) + B_0,
\end{equation}
where $B, B_0$ are some positive real numbers.

It can be assumed without loss of generality  that
$\vert {\mu}_k \vert > \max (1, 4d_2)$
for all
$k \in {\mathbb N}$ and
$\vert {\lambda}_j \vert > \max(2(d_1 +d_2),1)$ for all $j \in {\mathbb N}$.

At first, having functions  $L$ and $Y$, we construct
an entire function ${\cal N} \in P$, which outside
some  exceptional disks satisfies the estimate
$$
\vert {\cal N}(z) \vert \ge
C_{\cal N}
\exp(\psi( Im \ z ) + w_{m+1}(\vert z \vert)),
$$
where $ C_{\cal N}$ is some positive number.

By the maximum principle applied to $Y$ in disks
$D({\mu}_k, d_2 ),  k \in {\mathbb N}$,  and  using inequality (4),
we get from  (28)
$$
\vert Y(z) \vert
\le B_1 \exp(
 w_{m+3} (\vert \ z \vert) + B \ln(1+ \vert z \vert )), \ z \in
{\mathbb C},
$$
where $B_1$ is some positive number.
Whence,  using Lemma 2, we have for some $B_2 > 0$
\begin{equation}
\vert Y(z) \vert
\le B_2 \exp(
 w_{m+4} (\vert \ z \vert)), \ z \in {\mathbb C}.
\end{equation}

For any $\nu >0, b \in (0, d_2)$ let
$\Omega(b, \nu) = \bigcup \limits_{k=1}^{\infty} D({\mu}_k,
b {\vert {\mu}_k \vert}^{-\nu})$.
Let us estimate $\vert Y(z) \vert $ from below   outside the set
$\Omega(b, \nu) $. For any $k \in {\mathbb N}$
$$
Y_k(z)= \displaystyle \frac {z - {\mu}_k} {Y(z)} \ , \ z \in
\overline{D({\mu}_k, d_2)}.
$$
 $Y_k(z)$ is holomorphic in   $\overline
{D({\mu}_k, d_2)}$  and has no zeros there.

By the maximum principle
$$
\vert Y_k(z) \vert \le
\displaystyle \frac {d_2}{
\vert Y(z_o) \vert}
\ , \ z \in D({\mu}_k, d_2),
$$
where $z_0$ is a point of $ \in  \partial  D({\mu}_k,d_2)$
where  the minimum  of
$ \vert Y(z) \vert $ over
$ \overline {D({\mu}_k,d_2)} $ is attained.
>From this,  using  (28) and  (4), we obtain
$$
\vert Y_k(z) \vert \le C_1 (1+\vert z \vert)^B \exp(-w_{m+3}(\vert z \vert)),
\ z \in D({\mu}_k, d_2),
$$
where $C_1 >0$ doesn't depend on  $k$.
Hence
$$
\vert Y(z) \vert \ge \displaystyle \frac {\vert z- {\mu}_k \vert
\exp(w_{m+3}(\vert z \vert))}{C_1 (1+\vert z \vert)^B} \ ,
\ z \in D({\mu}_k, d_2).
$$
Thus, for
$
z \in D({\mu}_k, d_2)
\setminus D({\mu}_k, b {\vert {\mu}_k \vert}^{-\nu})
$
$$
\vert Y(z) \vert \ge \displaystyle \frac {b
\exp(w_{m+3}(\vert z \vert))}
{C_1 {\vert {\mu}_k \vert}^{\nu}(1+\vert z \vert)^B}.
$$
Since for
$ z \in D({\mu}_k, d_2)$ \ $\vert {\mu}_k \vert \le
\ \displaystyle \frac {4 \vert z \vert }{3}$,
then
\begin{equation}
\vert Y(z) \vert \ge \displaystyle \frac {C_2
\exp(w_{m+3}(\vert z \vert))}{(1+\vert z \vert)^{B+ \nu}}, \
z \in D({\mu}_k, d_2) \setminus D({\mu}_k,
 b {\vert{\mu}_k \vert}^{-\nu}) \ ,
\end{equation}
where $C_2 >0 $  doesn't depend on  $k$.

>From  (28) and  (30) we get
$$
\vert Y(z) \vert \ge \displaystyle \frac {C_3
\exp(w_{m+3}(\vert z \vert))}{(1+\vert z \vert)^{B+ \nu}}, \
z \notin
\Omega(b, \nu),
$$
where $C_3 >0$ is some constant. Using  Lemma 2, we get
\begin{equation}
\vert Y(z) \vert \ge C_4
\exp(w_{m+2}(\vert z \vert)), \ z \notin \Omega(b, \nu),
\end{equation}
where $C_4 >0 $ is some constant depending on $b, \nu $.

Now we need  to choose positive  numbers  $b,  \nu$
and sequence  $ t= \{t_j \}_{j=1}^{\infty} \subset {\mathbb C}$  such that:

1). for every $
j \in {\mathbb N}$ \ $D({\lambda}_j +t_j,
b {\vert {\lambda}_j +t_j \vert}^{-\nu})
\subset
D({\lambda}_j, d_1 {\vert {\lambda}_j \vert }^{1-\alpha});$

2). $ \Omega(b, \nu) \cap
\left( \bigcup \limits_{j=1}^{\infty}
D({\lambda}_j +t_j,
b {\vert {\lambda}_j +t_j \vert}^{-\nu}) \right)
= \varnothing;
$

3). $\sum \limits_{j=1}^{\infty} \vert t_j \vert {\vert {\lambda}_j
\vert}^{\alpha -1} < \infty. $

We  need only  slightly to develop  the scheme which was introduced in
[13] (see also  [12, Chapter 2, $ \S 4$ ]).

Let $n_{\mu}(r)$ is the number of zeros of $ Y $ of modulos at most
$r, r >0.$ From the estimates  (29), (1) it follows that $Y$ is
an entire function of order 1 and finite type. Consequently, for some
$A_{\mu}>1$ \  $ n_{\mu}(r) < A_{\mu}r , \ r>0 $ [2, \$ 5].
Choose
$ \nu$ and  $b$ so that  $ \nu >2 \alpha +1$,
$ 0 < b < \min(d_2,  d_1 2^{-(\nu +6)} A_{\mu}^{-2})$.

Let us introduce the  notations:
$ \Delta(z) = \vert z \vert, z \in {\mathbb C},
\ D_{k,Y}
=
D({\mu}_k, b {\vert {\mu}_k \vert}^{-\nu}), \
D_{j, L}=
D({\lambda}_j, d_1 {\vert {\lambda}_j \vert}^{1-\alpha}), \ k, j \in {\mathbb
N}$.

Let $j \in {\mathbb N}$ be arbitrary but fixed. For
$k \in {\mathbb N}$ such that
$
\Delta (D_{k,Y})
\cap
\Delta
(D_{j, L})
\ne \varnothing
$
we have
$
\vert \vert {\mu}_k \vert - \vert {\lambda}_j \vert \vert <l $ ,
where  $l=d_1 +d_2$. Sum  $s$ of lengths of all intervals
$\Delta (D_{k,Y})$  intersecting  $\Delta(D_{j, L})$
is estimated as follows:
$$s = \sum \limits_{\vert {\lambda}_j \vert - l <
\vert {\mu}_k \vert < \vert {\lambda}_j \vert
+l} 2b {\vert {\mu}_k \vert}^{-\nu}
\le \displaystyle \frac {2bn_{\mu}(\vert {\lambda}_j \vert +l)}{(\vert
{\lambda}_j \vert - l)^{\nu}} \le
\displaystyle \frac {2bA_{\mu}(\vert {\lambda}_j \vert +l)}
{(\vert {\lambda}_j \vert - l)^{\nu}}.
$$
Since
$\vert {\lambda}_j \vert > 2l$ for all $j \in {\mathbb N}$, then  $s\le
2^{\nu +2}
b A_{\mu} {\vert {\lambda}_j \vert}^{1-\nu }$. Put
${\sigma}_j =
2^{\nu +2}b A_{\mu}{\vert {\lambda}_j \vert}^{1 -\nu }$.
There exists an interval  $(\vert {\lambda}_j \vert +k_0
{\sigma}_j,
\vert {\lambda}_j \vert +(k_0+1){\sigma}_j)$, where  $ k_0 +1 \le
4n
(\vert
{\lambda}_j
\vert +l)$, not intersecting
$
\bigcup \limits_{k=1}^{\infty}
\Delta
(D_{k,Y})$. Put
$t_j=
(k_0+ \displaystyle \frac 1 2){\sigma}_j
\exp(i \arg{\lambda}_j)$. Thus, the choice of
the points $t_j$ is the same as in [12], [13].
Now we estimate $ \vert t_j \vert $:
\begin{equation}
\vert t_j
\vert \le (k_0+1){\sigma}_j \le 4n(\vert{\lambda}_j \vert +l)
2^{\nu +2}b A_{\mu}{\vert {\lambda}_j \vert}^{1 -\nu } \le
2^{\nu +5} b A_{\mu}^2
{\vert {\lambda}_j \vert}^{2 -\nu }.
\end{equation}
Consequently,
$$
\sum \limits_{j=1}^{\infty} \vert t_j \vert
{\vert
{\lambda}_j \vert}^{\alpha -1 }
\le
2^{\nu +5} b A_{\mu}^2
\sum \limits_{j=1}^{\infty} {\vert
{\lambda}_j \vert}^{\alpha +1 - \nu } < \infty,
$$
since   \  $\nu -\alpha -1 > \alpha$ and the covergence exponent
of the zeros of $ L $
doesn't exceed  $\alpha$.
This means that the
condition 3) holds.

Due to (32), choice of  $\nu$, $b$ and since   $\vert {\lambda}_j
\vert >1 $ for all $j \in {\mathbb N}$,  the condition
1) holds too.

By the construction disk  $D({\lambda}_j +t_j, \displaystyle
\frac {{\sigma}_j} {2})$
 doesn't intersect
 $ \Omega(b,\nu)$. Since
$\displaystyle \frac {{\sigma}_j} 2 =  2^{\nu+1} b A_{\mu}
{\vert {\lambda}_j \vert}^{1-\nu } >
b (\vert {\lambda}_j \vert +\vert t_j \vert)^{-\nu} =
b{\vert
{\lambda}_j
+t_j \vert}^{-\nu},
$
then disk
$D({\lambda}_j + t_j, \linebreak
b {\vert {\lambda}_j
+t_j \vert}^{-\nu})$  doesn't intersect $ \Omega(b, \nu) $ too.
Thus, condition 2) is satisfied.

Denote
$ D({\lambda}_j +t_j, b {\vert {\lambda}_j +t_j \vert}^{-\nu})$
by
$D_{j,L}'\ ,j \in {\mathbb N}.$

According to theorem  S \  $L_t$ satisfies
$$
\vert \ln \vert L_t(z) \vert - \ln \vert L(z) \vert \vert \le C, \
z \notin
\bigcup \limits_{j=1}^{\infty} D_{j,L} \ ,
$$
where $C>0$ is some constant.
>From this and (27) it follows that for
$ z \notin
\bigcup \limits_{j=1}^{\infty} D_{j,L}
$
\begin{equation}
\vert \ln \vert L_t(z) \vert - \psi ( Im \ z ) \vert \le
A \ln(1 + \vert z \vert) + A_0 + C.
\end{equation}
Note that  from  (2) for  $ \ z_1, z_2 \in D_{j, L} $
\begin{equation}
\vert\psi (Im \ z_1) -     \psi (Im \ z_2) \vert \le A_1,
\end{equation}
where $ A_1 > 0 $ is some constant independent of
$j \in {\mathbb N} $.

Applying the maximum principle to  $ L_t $  in every disk
$D_{j,L}$ and  using (33), (34),
we obtain
\begin{equation}
\vert L_t(z) \vert \le A_2 (1 + \vert z \vert)^{A}
\exp(\psi(  Im \ z )), \ z \in
\bigcup \limits_{j=1}^{\infty} D_{j,L},
\end{equation}
where $A_2 > \exp(A_0 + C)$ is some  constant.

>From  (33) and  (35) it follows that (35) holds
for all $z \in {\mathbb C}$.
Using the inequalities  (27), (33), (34),   it is not hard
to show ( see, for examaple, [8, Chapter 2, \S 4])
that outside
$ \bigcup \limits_{j=1}^{\infty}
D_{j,L}' $
\begin{equation}
\vert L_t(z) \vert \ge \exp(\psi( Im \ z ) -A_3 \ln(1 +\vert z \vert)
- A_4),
\end{equation}
where
$A_3, A_4 $ are some positive numbers.
Thus,  the necessary  estimates of  $ \vert L_t \vert $ from above
in ${\mathbb C}$
and from below outside
$ \bigcup \limits_{j=1}^{\infty} D_{j,L}' $
are established.

By the choice of  $b, \nu $ the sums of the radii of the disjoint exceptional
disks $D_{j,L}' $ and   $D_{k,Y}$
are  finite. Hence there exists a  sequence
$\{l_n\}_{n=1}^{\infty}$  of circles                    e
$l_n
=\{z \in {\mathbb
C}: \vert z \vert = R_n \}$, $R_n \to \infty $ as $n \to \infty $,
which doesn't
intersect the disks $D_{j,L}' $ and   $D_{k,Y}$.

Put
$$
{\cal L}(z)= \displaystyle \frac {L_t(z)}
{ \prod \limits_{j=1}^{N}
(z - ({\lambda}_j  +t_j))}, \ z \in {\mathbb C},
$$
where $N  \in {\mathbb N}$ is such that  $N > A +2 $.
Then from (35), (36) for some  positive constants $A_5, A_6$
\begin{equation}
\vert {\cal L}(z) \vert \le \displaystyle \frac
{A_5 \exp (\psi ( Im \ z ))}
{1 + {\vert z \vert}^2},\ z \in {\mathbb C};
\end{equation}

\begin{equation}
\vert {\cal L}(z) \vert \ge \displaystyle \frac {A_6
\exp (\psi ( Im \ z ))}
{(1 + \vert z \vert)^{A_3 + N}}, \ z \notin \bigcup \limits_{j=N+1}^{\infty}
D_{j,L}'.
\end{equation}

Let ${\cal N}(z) = {\cal L}(z) Y(z), \ z \in {\mathbb C}$.
Let $\{a_n\}_{n=1}^{\infty}$ be the zeros of  $\cal N$ ordered by
non-decreasing modulus.
By   (31), (38) and Lemma 2
for
$ z $
outside   \linebreak
the set $\bigcup \limits_{n=1}^{\infty}
D(a_n, b {\vert a_n \vert}^{-\nu})$
\begin{equation}
\vert {\cal N}(z) \vert \ge
C_{\cal  N}
\exp(\psi( Im \ z ) + w_{m+1}(\vert z  \vert)),
\end{equation}
where  $ C_{\cal  N} $ is some positive number.
In paricular, this estimate holds on the circles
$l_n, n \in {\mathbb N} $.

Let $K$ be an arbitrary compact in  ${\mathbb C}$ \ and   \
$k_0 \in {\mathbb N}$ be  such that
 ${\overline K} \subset D(0,R_{k_0}) $.
For  $k \ge k_0$ let
$$
I(z) = \int_{l_k} \displaystyle \frac {F(\xi)}
{{\cal N}(\xi) (\xi -z)} \ d  \xi,
\ z \in   K.
$$
Then
\begin{equation}
I(z)= 2\pi i \sum \limits_
{\vert a_n \vert < R_k} \displaystyle \frac {F(a_n)}
{{\cal N}'(a_n)(a_n - z)} + 2 \pi i
\displaystyle \frac {F(z)}{{\cal N}(z)}, \ z \in K.
\end{equation}
Using  (26) and  (39) for  $ z \in K $ we have
$$
\vert I(z) \vert \le \int_{l_k}
\displaystyle \frac {c
\exp(\psi( Im \ \xi ) +w_m(\vert \xi \vert))}
{C_{\cal  N}  \ dist(K, l_k)
\exp(\psi( Im \ \xi )
+w_{m+1}(\vert \xi \vert))} \vert
d \xi \vert.
$$

Letting  $k \to \infty $ in the right-hand side of this
inequality and taking Lemma 2 into account, we get  $I(z) =0, \ z \in K$.

Letting  $ k \to \infty $ in (40), we obtain
\begin{equation}
F(z)=
\lim \limits_{k \to \infty} \sum \limits_
{\vert a_n \vert < R_k}
\displaystyle \frac {F(a_n)}
{{\cal N}'(a_n)}
\displaystyle \frac {{\cal N}(z)}
{z-a_n} \ , \ z \in K.
\end{equation}
Since compact  $K$ was arbitrary, then  (41) holds  for all
$
z \in {\mathbb C}
$.

Consider the series
\begin{equation}
\sum \limits_{n=1}^{\infty}
\displaystyle \frac {F(a_n)}
{{\cal N}'(a_n)}
\displaystyle \frac {{\cal N}(z)}
{z-a_n} \ , \ z \in {\mathbb C}.
\end{equation}
We wish to show that it converges uniformly on  every  compact set
of complex plane.
 At first  for any  $n \in {\mathbb N}$ we estimate
$ \displaystyle \frac {{\cal N}(z)}
{z-a_n},
\ z \in {\mathbb C}$.
We consider two cases.

The first case.
Let $a_n$ be the zero of ${\cal L} $. Then for some
$j \in {\mathbb N} $ \  $ a_n = {\lambda}_j + t_j$, hence
$$
\displaystyle \frac {{\cal N}(z)}{z-a_n}=
\displaystyle \frac {{\cal L}(z)}
{z-({\lambda}_j +t_j)} \ Y(z), \ z \in {\mathbb C}.
$$
If  $ z \notin D_{j,L}' $,
then by  (37)
\begin{equation}
\left\vert \displaystyle \frac {{\cal L}(z)}
{z-({\lambda}_j +t_j)} \right\vert \le
\displaystyle \frac
{A_5 \exp(\psi( Im \ z )) }{1 + {\vert z \vert}^2}
 \displaystyle \frac {{\vert {\lambda}_j + t_j \vert}^{\nu}} b
\ .
\end{equation}

If  $ z \in D_{j,L}'$,
then by the maximum principle
$$
\left\vert \displaystyle \frac {{\cal L}(z)}
{z-({\lambda}_j +t_j)} \right\vert \le
\max \limits_
{\xi \in \partial D_{j,L}'}
\left\vert \displaystyle \frac {{\cal L}(\xi)}
{\xi-({\lambda}_j +t_j)} \right\vert \ .
$$

Let maximum in the right-hand side of the last inequality is attained
at a point
${\xi}_0 \in \partial D_{j,L}'$. Then, again using  the inequality (37),
$$
\left\vert
\displaystyle \frac
{ {\cal L}(z)}
{z-({\lambda}_j +t_j)} \right\vert \le
\displaystyle \frac
{\vert {\cal L}({\xi}_0) \vert {\vert {\lambda}_j  +
t_j \vert}^{\nu}}
b
\le
 \displaystyle \frac
{A_5 \exp(\psi( Im \ {\xi}_0) ) }{1 + {\vert {\xi}_0 \vert}^2}
\displaystyle \frac
{{\vert {\lambda}_j  +  t_j \vert}^{\nu}} b \ .
$$

Taking into account that for  $ z \in D_{j,L}'$ \
$\vert {\xi}_0  \vert \ge \vert z \vert -  1  $,
and using (34), we obtain
\begin{equation}
\left\vert
\displaystyle \frac
{{\cal L}(z)}
{z-({\lambda}_j +t_j)} \right\vert \le
\displaystyle \frac
{A_7 \exp(\psi( Im \ z )) }{1 + {\vert z \vert}^2}
\displaystyle \frac
{{\vert {\lambda}_j  + t_j \vert}^{\nu}} b,
\end{equation}
where $A_7 >0$ is some constant independent of  $j$.
Thus, from   (43), (44) we have
\begin{equation}
\left\vert
\displaystyle \frac
{{\cal L}(z)}{z-a_n}
\right\vert \le
\displaystyle \frac {A_8
\exp(\psi( Im \ z )) }{1 + {\vert z \vert}^2}
{\vert a_n \vert}^{\nu},
\end{equation}
where $ A_8 = b^{-1} \max (A_5, A_7)$.
>From  (29), (45) we get in the first case

\begin{equation}
\left\vert
\displaystyle \frac
{{\cal N}(z)}{z-a_n} \right\vert
\le
A_9 {\vert a_n \vert}^{\nu}
\exp(\psi( Im \ z ) + w_{m+4}(\vert z \vert ) ),
\ z \in {\mathbb C},
\end{equation}
where  $A_9 = A_8 B_2 $.

The second case. Let $a_n $ be the zero of  $ Y $. Then for some
$k \in {\mathbb N} $ \ $ a_n = {\mu}_k $. Hence
$$
\displaystyle \frac {{\cal N}(z)}{z-a_n}=
\displaystyle \frac {Y(z)}
{z- {\mu}_k} \
{\cal L}(z), \ z \in {\mathbb C}.
$$

If  $ z \notin
D_{k,Y} $, then from  (29)
\begin{equation}
\left\vert \displaystyle \frac {Y(z)}
{z- {\mu}_k} \right\vert \le
b^{-1} B_2
\exp (w_{m+4} (\vert \ z \vert))
{\vert {\mu}_k \vert}^{\nu},
\ z \in {\mathbb C}.
\end{equation}
If $ z \in
D_{k,Y}
$, then by the
the maximum principle

$$
\left\vert \displaystyle \frac {Y(z)}
{z- {\mu}_k} \right\vert \le \max \limits_{\xi \in
\partial D_{k,Y}}
\left\vert \displaystyle \frac {Y(\xi) }
{\xi- {\mu}_k} \right\vert \ .
$$

Let maximum in the right-hand of the last inequality is attained
at a point
${\xi}_0 \in
\partial D_{k,Y}$. Then,  again using  (29), we have
$$
\left\vert
\displaystyle \frac {Y(z)}
{z- {\mu}_k} \right\vert \le
b^{-1}
\vert Y({\xi}_0) \vert
{\vert {\mu}_k \vert}^{\nu} \le
b^{-1}B_2
\exp(w_{m+4}(\vert {\xi}_0 \vert)
{\vert {\mu}_k \vert}^{\nu}  \ .
$$
For  $ z \in D_{k,Y}$ \  $\vert {\xi}_0 \vert \le
\vert z \vert  + 1 $, so from the last inequality,
using (4), we get
\begin{equation}
\left\vert
\displaystyle \frac {Y(z)}
{z- {\mu}_k} \right\vert \le
B_3\exp(w_{m+4}(\vert z \vert )
{\vert {\mu}_k \vert}^{\nu},
\end{equation}
where $B_3 >0$ doesn't depend on  $k$.
>From  (47), (48) it follows that
\begin{equation}
\left\vert
\displaystyle \frac {Y(z)}
{z- a_n} \right\vert
\le
B_4 {\vert a_n \vert}^{\nu}
\exp(w_{m+4}(\vert z \vert )),
\end{equation}
where $B_4= \max ( b^{-1} B_2, B_3) $ doesn't depend on $n$. From  (37)
and (49) in the second case we obtain
\begin{equation}
\left\vert
\displaystyle \frac {{\cal N}(z)}
{z- a_n} \right\vert \le
A_{10}
{\vert a_n \vert}^{\nu} \exp(\psi( Im \ z ) +w_{m+4}(\vert z \vert )),
\ z \in {\mathbb C},
\end{equation}
where  $A_{10}= A_5 B_4$.
>From  (46) and  (50) we conclude that in both cases
\begin{equation}
\left\vert
\displaystyle \frac {{\cal N}(z)}
{z- a_n} \right\vert \le
A_{11}
{\vert a_n \vert}^{\nu}
\exp(\psi( Im \ z )
+w_{m+4}(\vert z \vert )), \ z \in {\mathbb C},
\end{equation}
where $A_{11}= \max (A_9, A_{10}) $.

>From the representation
$$
\displaystyle \frac
1 {{\cal N}'(a_n)} = \displaystyle \frac
{1} {2 \pi i}
\int \limits_{\vert \xi - a_n \vert =
b {\vert a_n \vert}^{-\nu}}
\displaystyle \frac {d \ {\xi}}
{{\cal N}(\xi) } \ ,
$$
the inequality (39) and by using  (4) and  (34),  we get
$$
\left\vert
\displaystyle \frac
1 {{\cal N}'(a_n)} \right\vert \le
A_{12}{\vert a_n \vert}^{-\nu}
\exp(-(\psi (Im \ a_n ) + w_{m+1}(\vert a_n \vert ))),
$$
where constant $A_{12} >0$ doesn't depend on  $n \in {\mathbb N}$.
>From this
and (26)
we obtain
$$
\left\vert
\displaystyle \frac
 {F(a_n)} {{\cal N}'(a_n)} \right\vert \le
c A_{12}{\vert a_n \vert}^{-\nu}
\exp(w_m(\vert a_n \vert ) -w_{m+1}(\vert a_n \vert )), \ n \in {\mathbb N}.
$$
By Lemma 2 there exists constant
 $Q \ge 0 $ such that
$w_m(\vert a_n \vert ) -w_{m+1}(\vert a_n \vert ) < -(\alpha +1)
\ln (1 + \vert a_n \vert) + Q  \ \ \forall n \in {\mathbb N}.
$
Consequently,
\begin{equation}
\left\vert
\displaystyle \frac
{F(a_n)} {{\cal N}'(a_n)} \right\vert \le
A_{13}
{(1 + \vert a_n \vert)}^{-(\alpha + \nu + 1)}, \ n \in {\mathbb N},
\end{equation}
where constant $ A_{13}>0 $ doesn't depend on  $n$.

>From  (51) and  (52) for every $z \in {\mathbb C}$ we have
 $$
\sum \limits_{n=1}^{\infty}
\left\vert
\displaystyle \frac {F(a_n)}
{{\cal N}'(a_n)}
\right\vert
\left\vert
\displaystyle \frac {{\cal N}(z)}
{z-a_n}
\right\vert \le
A_{11} A_{13}
\sum \limits_{n=1}^{\infty}
{(1 + \vert a_n \vert)}^{-(\alpha  + 1)}
\exp (\psi( Im \ z ) +
w_{m+4}(\vert z \vert)).
$$
Since
the convergence exponent
of the zeros of  $ {\cal N}$ doesn't exceed  ${\alpha}$, then
\begin{equation}
\sum \limits_{n=1}^{\infty}
{(1 + \vert a_n \vert)}^{-(\alpha  + 1)} < \infty
\end{equation}
Consequently, the series  (42) converges uniformly on every compact subset
of   ${\mathbb C}$.
Hence by (41) ((41) holds for all $ z \in {\mathbb C} $) we have
\begin{equation}
F(z) =
\displaystyle \sum \limits_{n=1}^{\infty}
\displaystyle \frac {F(a_n)}
{{\cal N}'(a_n)}
\displaystyle \frac {{\cal N}(z)}
{z-a_n}
\ ,
z \in {\mathbb C}.
\end{equation}

Now we wish to define
${\cal F}_n \in G^*$
such that
$ {\hat {\cal F}_n}(z) =
\displaystyle \frac {{\cal N}(z)}{z-a_n} \ , \linebreak
z \in {\mathbb C}, n \in {\mathbb N} $.
There are two cases.

The first case. $a_n$ is the zero of  $
{\cal L}
$. Hence,$
\displaystyle \frac {{\cal N}(z)}{z-a_n} =
\displaystyle \frac {{\cal L}(z)}{z-a_n} Y(z), z \in {\mathbb C}.
$
By (29), (45), Lemma 6
there exists
$
{\cal F}_n
\in G^*$ such that
$ {\hat {\cal F}_n}(z) = \displaystyle \frac {{\cal N}(z)}{z-a_n}
\ ,
z \in
{\mathbb C}$ and
$
\vert {\cal F}_n (f)
 \vert \le {\beta}_{m+4} A_{8} B_2
{\vert a_n \vert}^{\nu}
{\Vert f \Vert }_{G,m+5}\ , \ f \in G.
$

The second case.  $a_n$ is the zero of $Y$. Hence,
$
\displaystyle \frac
{{\cal N}(z)}{z-a_n} =
\displaystyle \frac {Y(z)}{z-a_n}{\cal L}(z)\ , z \in {\mathbb C}.
$
By  (37), (49), Lemma 6
there is
$
{\cal F}_n
\in G^*$ such that
$ {\hat {\cal F}_n}(z) = \displaystyle \frac {{\cal N}(z)}{z-a_n}
\ ,
z \in
{\mathbb C}$
and
$
\vert {\cal F}_n (f)
 \vert \le {\beta}_{m+4} A_5 B_4
{\vert a_n \vert}^{\nu}
{\Vert f \Vert }_{G,m+5}\ , \ f \in G.
$

Thus, in both cases  for every $ n \in {\mathbb N}$
there is
${\cal F}_n \in G^*$
such that
\begin{equation}
\hat{{\cal F}_n}(z)=
\displaystyle \frac {{\cal N}(z)}{z-a_n} \ , z \in {\mathbb C},
\end{equation}
and
\begin{equation}
\vert {\cal F}_n (f)\vert \le H
{\vert a_n \vert}^{\nu}
{\Vert f \Vert}_{G,m+5} \ , f \in G,
\end{equation}
where constant $H > 0 $  doesn't depend on  $n$.
Put
$$
{\cal F}(f) =
\sum \limits_{n=1}^{\infty}\displaystyle \frac {F(a_n)}{{\cal N}'(a_n)}
{\cal F}_n (f), \ f \in G.
$$
Because of  (52), (53), (56) ${\cal F}$ is defined correctly
and  it belongs to $G^*$.
>From (54) and  (55) we have
${\hat {\cal F}}= F$.  Thus, $J$ is surjective.

By the open mapping  theorem for the spaces
$(LN^*)$ ([9], [6, p. 12]) \ $J$ establishes
topological isomorphism of the spaces  $G^*$ and  $P$.  Theorem 1 is proved.

Acknowledgement. The work was supported by the Russian Foundation of
Fundamental Researches under Grant  99-01-00655.
I am sincerely grateful to V. V. Napalkov for his interest in the work.

\pagebreak

\begin{center}
REFERENCES
\end{center}

1. H. Komatsu, Ultradistributions I. Structure theorems and a

characterization. {\it J. Fac. Sci. Tokyo Sec. IA }, {\bf 20} (1973), 607 -- 628.

2. B. Ya. Levin, {\it Disribution of zeros of entire functions},
Gostekhizdat,

Moscow, 1956; English transl., American Math. Soc.,
Providence, 1964.

3.  S. Mandelbrojt, {\it S\'eries adh\'erentes. R\'egularisation des suites.

Applications}, Paris, 1952; Russian transl., Inostrannaya Literatura,

Moscow,
1955.

4. M. M. Mannanov,  Description of a class of analytic  functionals,

{\it Sibirsk. Mat. Zh.} {\bf 31}  (1990), 62-72.
[In Russian].

5. V. V. Napalkov,  Spaces of analytic functions of given
growth near

the boundary, {\it Izv. Akad. Nauk SSSR  Ser. Mat.} {\bf 51} (1987),
287-305.

[In Russian].

6. V.V. Napalkov.
{\it Convolution equations in multidimensional spaces},

"Nauka", Moscow, 1982.
[In Russian].

7. V. V. Napalkov and M. I. Solomesch, Estimate of changing of entire

function
under shift of its zeros, {\it Dokl. of Akad. Nauk (Russia)}, {\bf 342}

(1995),
739 -- 741.

8. S. V. Popyonov, On a weighted space of functions analytic in the

unbounded convex domain in  ${\mathbb C}^m $,
{\it Mat. Zametki} {\bf 40}  (1986), 374-384.

[In Russian].

9. A. P. Robertson and W. Robertson, {\it Topological vector spaces},

Cambridge University Press, Cambridge, 1964.

10. R. T. Rockafellar,  {\it  Convex analysis},
Princeton  Univ. Press, Princeton,

NY, 1970; Russian transl., Mir, Moscow, 1973.

11.  J. Sebasti\~ao e Silva,  Su certe classi di spazi localmente convessi

importante per le applicazioni, {\it  Rend. Mat. e Appl.}
{\bf 141} (1955), 388-410;

Russian transl., {\it Matematika} {\bf 1} (1957), 60-77.

12.  M. I. Solomesch,  {\it Convolution type operators on some spaces of

analytic functions}, Dissertation, Institute of Mathematics with

Computing Center, Ufa Scientific Center, Russian Academy of sciences,

Ufa, 1995. [In Russian].

13. M. I. Solomesch,
On the R. S. Yulmukhametov's theorem
on the

approximation
of subharmonic functions,   Manuscript No. 2447 -- B92,

deposited with VINITI,
 24.07.1992.
[In Russian].

14. B. A. Taylor,   Analytically uniform spaces of infinitely differentiable

functions, {\it  Communications on pure and applied mathematics}
{\bf 24}  (1971),

39-51.

15. R. S. Yulmukhametov,
Quasianalytic classes of functions in
convex

domains, {\it Mat. Sb.} {\bf 130(172)} (1986), 500-519.
[In Russian].

16. R. S. Yulmukhametov, Approximation of subharmonic functions,

{\it Analysis Mathematica} {\bf 11}  (1985), 257-282.
[In Russian].

17. R. S. Yulmukhametov, Approximation of subharmonic functions,

{\it Mat. Sb.} {\bf 124(166)}  (1984), 393-415.
[In Russian].

\vspace {0.5 cm}

Il'dar Kh. Musin

Institute of mathematics,  Russian Academy of Sciences,
Chernyshevskii str., 112, Ufa,  450000, RUSSIA.

E-mail: musin@imat.rb.ru

\end{document}